\documentclass[11pt]{amsart} 
\usepackage{amssymb,amsmath,latexsym,enumerate,graphicx,verbatim} 

\def\Z{\mathbb{Z}}
\def\R{\mathbb{R}} 
 
\def\S{R} 
\def\Q{Q} 
\def\Ql{Q_\text{low}} 
\def\Qu{Q_\text{up}} 
\newcommand\fl[1]{\left\lfloor {#1} \right\rfloor} 
\newcommand\st[1]{\left( \left( {#1} \right) \right)} 
\newcommand\sts[1]{\left( \left( {#1} \right) \right)^\star} 
\renewcommand\sigma{S} 

\newtheorem*{theorem}{Theorem} 
\newtheorem*{proposition}{Proposition} 
\newtheorem*{algorithm}{Algorithm} 

\title{Refined upper bounds for the linear Diophantine problem of Frobenius} 

\author{Matthias Beck}
\address{Department of Mathematical Sciences\\
        Binghamton University (SUNY)\\
        Binghamton, NY 13902-6000\\
        USA}
\email{matthias@math.binghamton.edu}

\author{Shelemyahu Zacks}
\address{Department of Mathematical Sciences\\
        Binghamton University (SUNY)\\
        Binghamton, NY 13902-6000\\
        USA}
\email{shelly@math.binghamton.edu}

\copyrightinfo{}{}
\keywords{The linear Diophantine problem of Frobenius, upper bounds, Dedekind-Rademacher sums, reciprocity laws} 
\subjclass[2000]{11D04, 
                 05A15, 
                 11Y16} 

\begin{document}
\setlength{\parindent}{0pt} 

\abstract We study the \emph{Frobenius problem}: given relatively prime positive integers 
$a_1,\dots,a_d$, find the largest value of $t$ (the \emph{Frobenius number} $g(a_1,\dots,a_d)$) 
such that $ \sum_{k=1}^d m_{k} a_{k} = t $ has no  solution in nonnegative integers 
$ m_{ 1 } , \dots , m_{ d } $. We introduce a method to compute upper bounds for $g(a_1,a_2,a_3)$, 
which seem to grow considerably slower than previously known bounds. 
Our computations are based on a formula for the restricted partition function, which involves 
Dedekind-Rademacher sums, and the reciprocity law for these sums. 
\endabstract 

\maketitle  

\setlength{\parskip}{0.4cm}
\bibliographystyle{amsplain}


\section{Introduction}

Given positive integers $ a_{1} < a_2 < \dots < a_{d} $ with $ \gcd (a_{1}, \dots, a_{d}) = 1 $, 
the \emph{linear Diophantine problem of Frobenius} asks for 
the largest integer $t$ for which we cannot find 
nonnegative integers $ m_{ 1 } , \dots , m_{ d } $ such that
  \[  t = m_1 a_1 + \dots + m_d a_d \ . \] 
We call this largest integer the \emph{Frobenius number} $ g ( a_{1} , \dots , a_{d} ) $; 
its study was initiated in the 19th century. 
One fact which makes this problem attractive is that it can be easily
described, for example, in terms of coins of denominations $a_1, \dots, a_d$; the 
Frobenius number is the largest amount of money which cannot be formed using these coins. 
For $d=2$, it is well known (most probably at least since Sylvester \cite{sylvester}) that 
  \begin{equation}\label{gfor2} g(a_{1}, a_{2}) = a_{1} a_{2} - a_1 - a_2 \ . \end{equation} 
For $d>2$, all attempts to find explicit formulas have proved elusive. 
Two excellent survey papers on the Frobenius problem are \cite{alfonsin} and \cite{selmer}. 

Our goal is to establish upper bounds for $ g ( a_{1} , \dots , a_{d} ) $. 
The literature on such bounds is vast; it includes results by 
Erd\H{o}s and Graham \cite{erdosgrahamfrob} 
  \begin{equation}\label{knownbound1} g(a_{1}, \dots, a_{d}) \leq 2 a_d \left\lfloor \frac {a_1} d \right\rfloor - a_1 \ , \end{equation} 
Selmer \cite{selmer} 
  \begin{equation}\label{knownbound2} g(a_{1}, \dots, a_{d}) \leq 2 a_{d-1} \left\lfloor \frac {a_d} d \right\rfloor - a_d \ , \end{equation} 
and Vitek \cite{vitek} 
  \begin{equation}\label{knownbound3} g(a_{1}, \dots, a_{d}) \leq \left\lfloor \frac 1 2 ( a_2 - 1 ) ( a_d - 2 ) \right\rfloor - 1 \ . \end{equation} 
Here $a_1 < a_2 < \dots < a_d$, and $\lfloor x \rfloor$ denotes the greatest integer 
not exceeding $x$. 
Davison \cite{davison} established the lower bound 
  \begin{equation}\label{lowerbound} g ( a_{1} , a_2 , a_{3} ) \geq \sqrt{ 3 a_1 a_2 a_3 } - a_1 - a_2 - a_3 \ . \end{equation} 
Experimental data \cite{frobcomp} shows that Davison's bound is sharp in the sense 
that it is very often very close to $ g ( a_{1} , a_2 , a_{3} ) $. 
On the other hand, the upper bounds given by (\ref{knownbound1}), (\ref{knownbound2}), and (\ref{knownbound3}) 
seem to be quite large compared to the actual 
Frobenius numbers. In this paper, we derive a method of achieving sharper upper 
bounds for the Frobenius number. 
Our results are based on a formula for the restricted partition function (Section 
\ref{partitionsection}), which involves Dedekind-Rademacher sums, 
and the reciprocity law for these sums (Section \ref{rademachersection}). 
The main result is derived in Section \ref{mainsection}; 
computations which illustrate our new bounds can be found in Section \ref{compsection}. 

We focus on the first non-trivial case $d=3$; any bound for this case yields a general bound, 
as one can easily see that $ g ( a_{1} , \dots , a_{d} ) \leq g ( a_{1} , a_2 , a_{3} ) $ 
if $a_1$, $a_2$, and $a_3$ are relatively prime. If not then we can reduce by 
one variable at a time: Again by the definition of the Frobenius number, 
$g(a_1, \dots, a_d) \leq g(a_1, \dots, a_{d-1})$ if $a_1, \dots, a_{d-1}$ are 
relatively prime. If not, we can use a formula of Brauer and Shockely \cite{brauershockley}: 
If $n = \gcd \left( a_1, \dots, a_{d-1} \right)$ then 
  \begin{equation}\label{brauershock} g(a_1,\dots,a_d) = n \ g \left( \frac{a_1}{n} , \dots , \frac{a_{d-1}}{n} , a_d \right) + \left( n-1 \right) a_d \ . \end{equation} 
Hence 
  \[ g(a_1,\dots,a_d) \leq n \ g \left( \frac{a_1}{n} , \dots , \frac{a_{d-1}}{n} \right) + \left( n-1 \right) a_d \ . \] 


\section{The restricted partition function}\label{partitionsection} 

We approach the Frobenius problem through the study of the \emph{restricted partition function}
  \[ p_{ \{ a_1, \dots, a_d \} } (n) = \# \left\{ ( m_1, \dots, m_d ) \in \Z_{\geq 0}^d : \ m_1 a_1 + \dots + m_d a_d = n \right\} \ , \] 
the number of partitions of $n$ using only $a_1, \dots, a_d$ as parts. 
In view of this function, the Frobenius number $ g ( a_1, \dots, a_d ) $ is the largest integer $n$ 
such that $p_{ \{ a_1, \dots, a_d \} } (n) = 0$. 

In the $d=3$ case, we can additionally assume that $a=a_1$, $b=a_2$, and $c=a_3$ 
are \emph{pairwise} relatively prime, a simplification due to 
Johnson's formula \cite{johnson}: if $ n = \gcd (a,b) $ then 
  \[ g(a,b,c) = n \ g \left( \frac{a}{n} , \frac{b}{n} , c \right) + (n-1) c \ . \] 
(This identity is a special case of (\ref{brauershock}).) 

In the case that $ a, b, c $ are pairwise relatively prime, Beck, Diaz, and Robins derived 
the following result for the  partition function $p_{ \{ a, b, c \} }$ \cite[Theorem 3]{bdr}: 
  \begin{align} 
    p_{ \{ a, b, c \} } (n) &= \frac{ n^{ 2 }  }{ 2abc } + \frac{ n }{ 2 } \left( \frac{ 1 }{ ab } + \frac{ 1 }{ ac } + \frac{ 1 }{ bc }  \right) + \frac{ 1 }{ 12 } \left( \frac{ a }{ bc } + \frac{ b }{ ac } + \frac{ c }{ ab } \right) \label{p_abc} \\ 
                            &\qquad + \sigma_{-n} ( b, c; a ) + \sigma_{-n} ( c, a; b ) + \sigma_{-n} ( a, b; c ) \ . \nonumber 
  \end{align} 
Here \cite[Equation (14)]{bdr}
  \begin{equation}\label{fouriertriple} \sigma_{t} ( a, b; c ) = \sum_{ m=0 }^{ c-1 } \st{ \frac{ - a^{-1} ( b m + t ) }{ c } } \st{ \frac{ m }{ c } } \ , \end{equation} 
where $ a a^{-1} \equiv 1 \mod c$ and $\st x = x - \fl x - 1/2$, 
is a special case of a \emph{Dedekind-Rademacher sum}; we will discuss 
these sums in the next section. 

To bound the Frobenius number (from above), we need to bound $p_{ \{ a, b, c \} }$ 
(from below), whose only nontrivial ingredients are the Dedekind-Rademacher sums. 
A classical bound for the Dedekind-Rademacher sum yielded in \cite{bdr} the inequality 
  \[ g(a,b,c) \leq \frac{ 1 }{ 2 } \left( \sqrt{ abc \left( a + b + c \right) } - a - b - c \right) \ , \] 
which is of comparable size to the other upper bounds given by (\ref{knownbound1}), (\ref{knownbound2}), and (\ref{knownbound3}). 
However, we will show that one can obtain bounds of smaller magnitude. 


\section{Dedekind-Rademacher sums}\label{rademachersection} 

The \emph{Dedekind-Rademacher sum} \cite{rademacherdedekind} is defined for 
$ a, b \in \Z , \, x, y \in \R $ as 
  \[ \S (a,b;x,y) = \sum_{ k=0 }^{ b-1 } \sts{ \frac{ a (k+y) }{ b } + x } \sts{ \frac{ k+y }{ b } } \ , \] 
where 
  \[ \sts x = \left\{ \begin{array}{ll} \st x & \mbox{ if } x \not\in \Z , \\ 
                                        0     & \mbox{ if } x \in \Z . \end{array} \right. \] 
Rademacher's sum generalizes the classical Dedekind sum $ \S (a,b;0,0) $ \cite{dedekind}. 
An easy bound for the Dedekind-Rademacher sum $ \S (a,b;x,0) $ can be obtained through the 
Cauchy-Schwartz inequality: if $a$ and $b$ are relatively prime then 
  \begin{align} \left| \S (a,b;x,0) \right| &= \left| \sum_{ k=0 }^{ b-1 } \sts{ \frac{ a k }{ b } + x } \sts{ \frac{ k }{ b } } \right| \nonumber \\ 
                                            &\leq \sqrt{ \left( \sum_{ k=0 }^{ b-1 } { \sts{ \frac{ a k }{ b } + x } }^2 \right) \left( \sum_{ k=0 }^{ b-1 } { \sts{ \frac{ k }{ b } } }^2 \right) } \nonumber \\ 
                                            &= \sqrt{ \left( \sum_{ k=0 }^{ b-1 } { \sts{ \frac{ k }{ b } + x } }^2 \right) \left( \sum_{ k=1 }^{ b-1 } \left( \frac{ k }{ b } - \frac 1 2 \right)^2 \right) } \label{rademacherbound} \\ 
                                            &\leq \sqrt{ \left( \sum_{ k=0 }^{ b-1 } \left( \frac{ k }{ b } + \frac 1 b - \frac 1 2 \right)^2 \right) \left( \frac{ b }{ 12 } - \frac 1 4 + \frac{ 1 }{ 6 b } \right) } \nonumber \\ 
                                            &= \sqrt{ \left( \frac{ b }{ 12 } + \frac{ 1 }{ 6 b } \right) \left( \frac{ b }{ 12 } - \frac 1 4 + \frac{ 1 }{ 6 b } \right) } \nonumber 
  \end{align} 
(In the third and fourth step we use the periodicity of $\sts{x}$.) 
An important property of $ \S (a,b;x,y) $ is Rademacher's \emph{reciprocity law} 
\cite{rademacherdedekind}: if $a$ and $b$ are relatively prime then 
  \begin{equation}\label{rademacherreciprocity} \S (a,b;x,y) + \S (b,a;y,x) = \Q (a,b;x,y) \ . \end{equation} 
Here 
  \[ \Q (a,b;x,y) = \left\{ \begin{array}{ll} - \frac{ 1 }{ 4 } + \frac{ 1 }{ 12 } \left( \frac{ a }{ b } + \frac{ 1 }{ ab } + \frac{ b }{ a }  \right) & \mbox{ if both } x,y \in \Z , \\ 
                                              \sts x \sts y + \frac{ 1 }{ 2 } \left( \frac{ a }{ b } \psi_{ 2 } (y) + \frac{ 1 }{ ab } \psi_{ 2 } (ay+bx) + \frac{ b }{ a } \psi_{ 2 } (x) \right) & \mbox{ otherwise, } \end{array} \right. \] 
where 
  \[ \psi_{ 2 } (x) = ( x - \fl x )^{ 2 } - ( x - \fl x ) + 1/6 \] 
denotes the periodic second Bernoulli function. 
Among other things, this reciprocity law allows us to compute $ \S (a,b;x,y) $ 
in polynomial time, by means of a Euclidean-type algorithm using the first two variables: 
simply note that we can replace $a$ in $ \S (a,b;x,y) $ by the least residue of $a$ modulo $b$. 

To express $\sigma$ in terms of $\S$, we rewrite (\ref{fouriertriple}) as 
  \begin{align*} \sigma_t (a,b;c) &= \sum_{ m=0 }^{ c-1 } \sts{ \frac{ - a^{-1} ( b m + t ) }{ c } } \sts{ \frac{ m }{ c } } + \left\{ \begin{array}{ll} \frac 1 4 & \text{ if } c|t , \\ 
                                                                      - \frac 1 2 \sts{ - \frac{ a^{-1} t }{ c } } - \frac 1 2 \sts{ - \frac{ b^{-1} t }{ c } } & \text{ otherwise. } \end{array} \right. \end{align*} 
Accordingly, 
  \[ \sigma_t (a,b;c) = \S \left( - a^{-1} b, c; - \frac{ a^{-1} t }{ c } , 0 \right) 
     + \left\{ \begin{array}{ll} \frac 1 4 & \text{ if } c|t , \\ 
                                 - \frac 1 2 \sts{ - \frac{ a^{-1} t }{ c } } - \frac 1 2 \sts{ - \frac{ b^{-1} t }{ c } } & \text{ otherwise. } \end{array} \right. \] 
To ease our computations, we bound this as 
  \begin{equation}\label{sigmatrivbound} \sigma_t (a,b;c) \geq \S \left( - a^{-1} b, c; - \frac{ a^{-1} t }{ c } , 0 \right) - \frac 1 2 \ . \end{equation} 


\section{Upper bounds for $g(a,b,c)$}\label{mainsection} 

To bound $\sigma_t (a,b;c)$ from below (which yields an upper bound for $g(a,b,c)$), we use 
an interplay of (\ref{rademacherreciprocity}) and (\ref{rademacherbound}) to obtain a bound 
for the Dedekind-Rademacher sum corresponding to $\sigma_t$, according to (\ref{sigmatrivbound}). 
The idea is to reduce the arguments of the Dedekind-Rademacher sum after the application 
of (\ref{rademacherreciprocity}), which means that the bound given by (\ref{rademacherbound}) 
will be more accurate. To illustrate this, let $c_1$ be the least nonnegative residue of 
$- a^{-1} b$ modulo $c$. Then 
  \begin{equation}\label{firstrefinement} \S \left( - a^{-1} b, c; - \frac{ a^{-1} t }{ c } , 0 \right) = \S \left( c_1 , c; - \frac{ a^{-1} t }{ c } , 0 \right) = \Q \left( c_1 , c ; - \frac{ a^{-1} t }{ c } , 0 \right) - \S \left( c , c_1 ; 0 , - \frac{ a^{-1} t }{ c } \right) \ . \end{equation} 
If $c_1 = 1$ then the right-hand side can be simplified, as 
$ \S \left( c , 1 ; 0 , - \frac{ a^{-1} t }{ c } \right) = 0 $. 
If $c_1 \not= 1$ then the Dedekind-Rademacher sum on the right-hand side of 
(\ref{firstrefinement}) can be bounded (via (\ref{rademacherbound})) 
sharper then the Dedekind-Rademacher sum on the left-hand side. 
In fact, by a repeated application of (\ref{rademacherreciprocity}), we can achieve bounds 
which are even better. To keep the computations somewhat simple, we apply 
(\ref{rademacherreciprocity}) once more and illustrate what this process yields in terms of 
lower bounds for $\sigma_t$. Let $c_2$ be the least nonnegative residue of $c$ modulo $c_1$. 
If $c_2 = 1$ then 
  \begin{align} \S \left( - a^{-1} b, c; - \frac{ a^{-1} t }{ c } , 0 \right) &= \Q \left( c_1 , c ; - \frac{ a^{-1} t }{ c } , 0 \right) - \S \left( c , c_1 ; 0 , - \frac{ a^{-1} t }{ c } \right) \nonumber \\ 
    &= \Q \left( c_1 , c ; - \frac{ a^{-1} t }{ c } , 0 \right) - \S \left( 1 , c_1 ; 0 , - \frac{ a^{-1} t }{ c } \right) \label{c2=1} \\ 
    &= \Q \left( c_1 , c ; - \frac{ a^{-1} t }{ c } , 0 \right) - \Q \left( 1 , c_1 ; 0 , - \frac{ a^{-1} t }{ c } \right) , \nonumber \end{align} 
as $\S \left( c_1 , 1 ; - \frac{ a^{-1} t }{ c } , 0 \right) = 0$. 
If $c_2 \not= 1$ then (\ref{firstrefinement}) can be refined as 
  \begin{align} \S \left( - a^{-1} b, c; - \frac{ a^{-1} t }{ c } , 0 \right) &= \Q \left( c_1 , c ; - \frac{ a^{-1} t }{ c } , 0 \right) - \S \left( c_2 , c_1 ; 0 , - \frac{ a^{-1} t }{ c } \right) \nonumber \\ 
    &= \Q \left( c_1 , c ; - \frac{ a^{-1} t }{ c } , 0 \right) - \Q \left( c_2 , c_1 ; 0 , - \frac{ a^{-1} t }{ c } \right) + \S \left( c_1 , c_2 ; - \frac{ a^{-1} t }{ c } , 0 \right) . \label{secondrefinement} \end{align} 
The Dedekind-Rademacher sum on the right-hand side can be bounded according to 
(\ref{rademacherbound}) as 
  \begin{equation}\label{rademacherbound2} \S \left( c_1 , c_2 ; - \frac{ a^{-1} t }{ c } , 0 \right) \geq - \sqrt{ \left( \frac{ c_2 }{ 12 } + \frac{ 1 }{ 6 c_2 } \right) \left( \frac{ c_2 }{ 12 } - \frac 1 4 + \frac{ 1 }{ 6 c_2 } \right) } \ . \end{equation}  
We still need to bound $\Q$. 
$\psi_{ 2 }$ has a minimum of $-1/12$ (at $x = 1/2$) and a maximum of $1/6$ (at $x=0$). 
These extreme values yield for 
  \[ \Q \left( c_1 , c ; - \frac{ a^{-1} t }{ c } , 0 \right) = \left\{ \begin{array}{ll} - \frac{ 1 }{ 4 } + \frac{ 1 }{ 12 } \left( \frac{ c_1 }{ c } + \frac{ 1 }{ c_1 c } + \frac{ c }{ c_1 } \right) & \mbox{ if } c|t , \\ 
                                                                \frac{ 1 }{ 2 } \left( \frac{ c_1 }{ 6c } + \frac{ 1 }{ 6 c_1 c } + \frac{ c }{ c_1 } \psi_{ 2 } \left( - \frac{ a^{-1} t }{ c } \right) \right) & \mbox{ otherwise, } \end{array} \right. \] 
the lower bound 
  \begin{equation}\label{Qlower} \Q \left( c_1 , c ; - \frac{ a^{-1} t }{ c } , 0 \right) \geq - \frac{ 1 }{ 4 } + \frac{ c_1 }{ 12 c } + \frac{ 1 }{ 12 c_1 c } - \frac{ c }{ 24 c_1 } = \Ql (c_1,c) \ , \end{equation}  
as well as the upper bound 
  \begin{equation}\label{Qupper} \Q \left( c_2 , c_1 ; 0, - \frac{ a^{-1} t }{ c } \right) \leq \frac{ c_2 }{ 12 c_1 } + \frac{ 1 }{ 12 c_2 c_1 } + \frac{ c_1 }{ 12 c_2 } = \Qu (c_2,c_1) \ . \end{equation} 
These inequalities yield the following. 
\begin{proposition} Suppose $a$ and $b$ are relatively prime to $c$. Let $c_1$ be 
the least nonnegative residue of $- a^{-1} b$ modulo $c$, and let $c_2$ be the least 
nonnegative residue of $c$ modulo $c_1$. 
\begin{enumerate}[{\rm (i)}] 
\item\label{i1} If $c_1 = 1$ then \ 
  $\displaystyle \sigma_t (a,b;c) \geq - \frac{ c }{ 24 } + \frac{1}{6c} - \frac 3 4 $ . 
\item\label{i2} If $c_1 \not= 1$ and $c_2 = 1$ then \ 
  $\displaystyle \sigma_t (a,b;c) \geq \frac{ c_1 }{ 12 c } + \frac{ 1 }{ 12 c_1 c } - \frac{ c }{ 24 c_1 } - \frac{ 1 }{ 6 c_1 } - \frac{ c_1 }{ 12 } - \frac 3 4 $ . 
\item\label{i3} If $c_1 \not= 1$ and $c_2 \not= 1$ then 
  \[ \sigma_t (a,b;c) \geq \frac{ c_1 }{ 12 c } + \frac{ 1 }{ 12 c_1 c } - \frac{ c }{ 24 c_1 } - \frac{ c_2 }{ 12 c_1 } - \frac{ 1 }{ 12 c_1 c_2 } - \frac{ c_1 }{ 12 c_2 } - \frac 3 4 - \sqrt{ \left( \frac{ c_2 }{ 12 } + \frac{ 1 }{ 6 c_2 } \right) \left( \frac{ c_2 }{ 12 } - \frac 1 4 + \frac{ 1 }{ 6 c_2 } \right) } \ . \] 
\end{enumerate} 
\end{proposition} 
\emph{Proof.} (\ref{i1}) Use (\ref{firstrefinement}) with $c_1=1$ in (\ref{sigmatrivbound}): 
  \begin{align*} \sigma_t (a,b;c) &\geq \S \left( - a^{-1} b, c; - \frac{ a^{-1} t }{ c } , 0 \right) - \frac 1 2 \\ 
                                  &= \Q \left( 1 , c ; - \frac{ a^{-1} t }{ c } , 0 \right) - \frac 1 2 \\ 
                                  &\geq - \frac{ 1 }{ 4 } + \frac{ 1 }{ 6 c } - \frac{ c }{ 24 } - \frac 1 2 \ . \end{align*} 
Here the last inequality follows from (\ref{Qlower}). 

(\ref{i2}) Use (\ref{c2=1}) in (\ref{sigmatrivbound}) together with the bounds (\ref{Qlower}) 
and (\ref{Qupper}): 
  \begin{align*} \sigma_t (a,b;c) &\geq \S \left( - a^{-1} b, c; - \frac{ a^{-1} t }{ c } , 0 \right) - \frac 1 2 \\ 
                                  &= \Q \left( c_1 , c ; - \frac{ a^{-1} t }{ c } , 0 \right) - \Q \left( 1 , c_1 ; 0 , - \frac{ a^{-1} t }{ c } \right) - \frac 1 2 \\ 
                                  &\geq - \frac{ 1 }{ 4 } + \frac{ c_1 }{ 12 c } + \frac{ 1 }{ 12 c_1 c } - \frac{ c }{ 24 c_1 } - \left( \frac{ 1 }{ 6 c_1 } + \frac{ c_1 }{ 12 } \right) - \frac 1 2 \end{align*} 
(\ref{i3}) Use (\ref{secondrefinement}) with the bounds given in (\ref{rademacherbound2}), 
(\ref{Qlower}), and (\ref{Qupper}): 
  \begin{align*} \sigma_t (a,b;c) &\geq \S \left( - a^{-1} b, c; - \frac{ a^{-1} t }{ c } , 0 \right) - \frac 1 2 \\ 
                                  &= \Q \left( c_1 , c ; - \frac{ a^{-1} t }{ c } , 0 \right) - \Q \left( c_2 , c_1 ; 0 , - \frac{ a^{-1} t }{ c } \right) + \S \left( c_1 , c_2 ; - \frac{ a^{-1} t }{ c } , 0 \right) - \frac 1 2 \\ 
                                  &\geq - \frac{ 1 }{ 4 } + \frac{ c_1 }{ 12 c } + \frac{ 1 }{ 12 c_1 c } - \frac{ c }{ 24 c_1 } - \left( \frac{ c_2 }{ 12 c_1 } + \frac{ 1 }{ 12 c_2 c_1 } + \frac{ c_1 }{ 12 c_2 } \right) \\ 
				  &\qquad - \sqrt{ \left( \frac{ c_2 }{ 12 } + \frac{ 1 }{ 6 c_2 } \right) \left( \frac{ c_2 }{ 12 } - \frac 1 4 + \frac{ 1 }{ 6 c_2 } \right) }  - \frac 1 2 \ . \qquad \Box \end{align*} 

These lower bounds can be combined with (\ref{p_abc}) and the quadratic formula to give an upper bound 
on the Frobenius number. 
\begin{theorem}\label{mainthm} Suppose $a$, $b$, and $c$ are pairwise relatively prime. Denote the lower 
bounds for $\sigma_t (b,c;a)$, $\sigma_t (c,a;b)$, and $\sigma_t (a,b;c)$ according to the 
previous proposition by $\alpha$, $\beta$, and $\gamma$, respectively. Then 
  \[ g(a,b,c) \ \leq \ \sqrt{ \tfrac 1 4 \left( a+b+c \right)^2 - \tfrac 1 6 \left( a^2 + b^2 + c^2 \right) - 2 abc \left( \alpha + \beta + \gamma \right) } \, - \, \tfrac 1 2 (a+b+c) . \] 
\end{theorem} 
One should note that $\alpha + \beta + \gamma$ is negative. We can see that the growth 
behavior of this upper bound is dominated by $- 2 abc \left( \alpha + \beta + \gamma \right)$ 
under the square root. This means that if we can make $- \left( \alpha + \beta + \gamma \right)$ somewhat smaller 
than $\min (a,b,c)$ then we get a bound which grows considerably less that the bounds given 
by (\ref{knownbound1}), (\ref{knownbound2}), and (\ref{knownbound3}). In fact, we can see 
this difference in example computations already when we use the bounds $\alpha, \beta, \gamma$ as given by our 
proposition. 
What is more important, however, is the fact that we can easily obtain even better bounds by improving 
our proposition through additional applications of Rademacher's reciprocity law (\ref{rademacherreciprocity}). 
We illustrate this with the following algorithm, whose result is a bound on $\sigma_t (a,b;c)$, 
which can be used in the above theorem (instead of the bounds coming from the proposition). 
\begin{algorithm} Input: $a,b,c$ (pairwise relatively prime) and $N$ (number of iterations). 
Output: lower bound $S$ for $\sigma_t (a,b;c)$. 
\begin{verbatim} 
c_1 := - a^{-1} b  modulo c (least nonnegative residue) 
S := 0 
n := 1 

REPEAT { 

  c_2 := c  modulo c_1 (least nonnegative residue) 
  S_1 := S + Q_low (c_1,c) 
  S_2 := S_1 - Q_up (c_2,c_1) 

  IF c_1 = 1 THEN S := S_1 
  Else S := S_2 

  IF c_1 = 1 OR c_2 = 1 OR n = N THEN BREAK 

  c := c_2
  c_1 := c_1  modulo c_2 (least nonnegative residue) 
  n := n + 1 

} 

IF c_1 > 1 AND c_2 > 2 THEN S := S - sqrt((c_2/12 + 1/(6 c_2) - 1/4) (c_2/12 + 1/(6 c_2))) 

S := S - 1/2 
\end{verbatim} 
\end{algorithm} 
The algorithm repeats the steps described in the proposition $N$ times, at each step bounding $Q$ 
coming from Rademacher reciprocity according to (\ref{Qlower}) and (\ref{Qupper}). It stops 
prematurely if one of the variables is 1, in which case the remaining Dedekind-Rademacher sum is 
zero.


\section{Computations}\label{compsection} 

In the present section we illustrate the newly proposed upper bound for $g(a,b,c)$ numerically. 
In order to compare the results also with the lower bound given by Davison (\ref{lowerbound}) 
we present here the values 
  \[ f(a,b,c) = g(a,b,c) + a + b + c \ . \] 
For these Frobenius numbers, Davison's lower bound is 
  \[ f(a,b,c) \geq \sqrt 3 \, z \ , \] 
where $z = \sqrt{abc}$. 
In \cite{frobcomp} we presented together with David Einstein an algorithm for the exact computation of $f(a,b,c)$. 
Einstein computed 20000 ``admissible" (see \cite{frobcomp}) values of $f(a,b,c)$ 
for relatively prime arguments chosen at random from the set $\{ 3, \dots, 750 \}$. 
In \cite{frobcomp} we arrived at the empirical conjecture that $ f(a,b,c) \leq \sqrt{abc}^{5/4} $. 
The objectives of our current presentations are: 
\begin{enumerate}[(i)] 
  \item to compare our new upper bound with the known upper bound, which is, according to (\ref{knownbound1}), (\ref{knownbound2}), and 
  (\ref{knownbound3}), 
  \[ \min \left( 2c \left\lfloor \frac a 3 \right\rfloor - a , 2b \left\lfloor \frac c 3 \right\rfloor - c , \left\lfloor \frac 1 2 ( b - 1 ) ( c - 2 ) \right\rfloor - 1 \right) + a + b + c \ ; \] 
  \item to compare our new upper bound with the conjectured upper bound $z^{5/4}$; 
  \item to compare the new upper bound to the true value of $f(a,b,c)$. 
\end{enumerate} 

\begin{figure}[htb]
\begin{center}
\includegraphics[totalheight=4in]{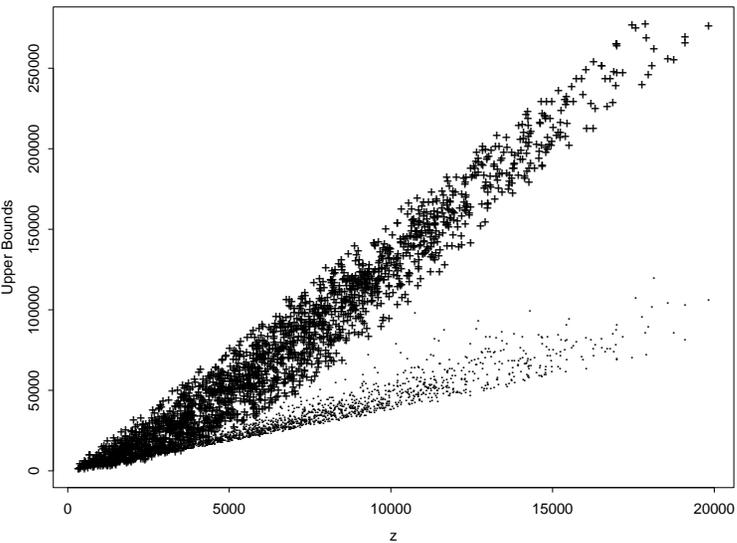}
\end{center}
\caption{The new and old upper bounds for the Frobenius number}\label{mainplot} 
\end{figure} 

For these objectives we computed the new upper bound and the known upper bound for two 
thousand values of $(a,b,c)$, randomly chosen from Einstein's data. 
In all computations we used the minimum of the lower bounds 
given by the proposition and the algorithm for $N=2$ in the theorem to obtain an upper 
bound for $f(a,b,c)$. 
In Figure \ref{mainplot} we plot the new upper bound ($\cdot$) 
and the known upper bound (+) as functions of $z$. 

\begin{figure}[htb]
\begin{center}
\includegraphics[totalheight=4in]{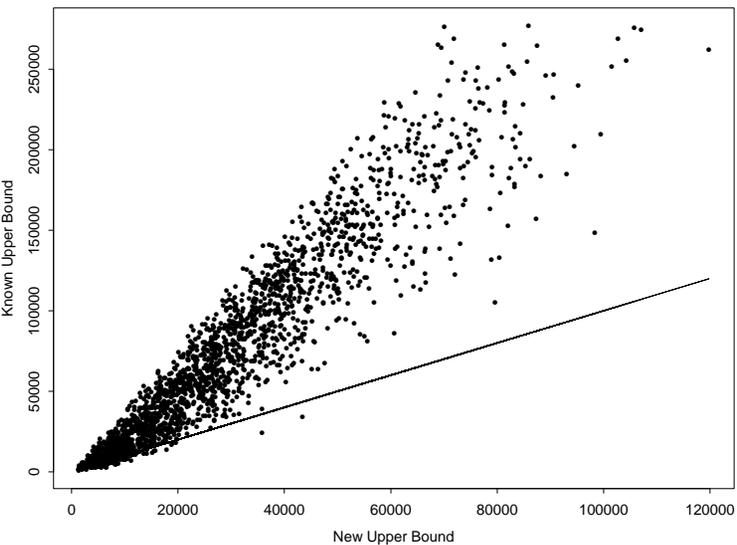}
\end{center}
\caption{The new and old upper bounds compared}\label{boundcomparisonfig} 
\end{figure} 

Among the 2000 points only less than 100 have known upper bounds smaller than the new 
upper bound. In 50\% of the cases the ratio of the known upper bound to the new upper bound 
is greater than 2.44. In Figure \ref{boundcomparisonfig} we plot the known upper bound 
as a function of the new upper bound. This figure complements the inferences from 
Figure \ref{mainplot}. 

\begin{figure}[htb]
\begin{center}
\includegraphics[totalheight=4in]{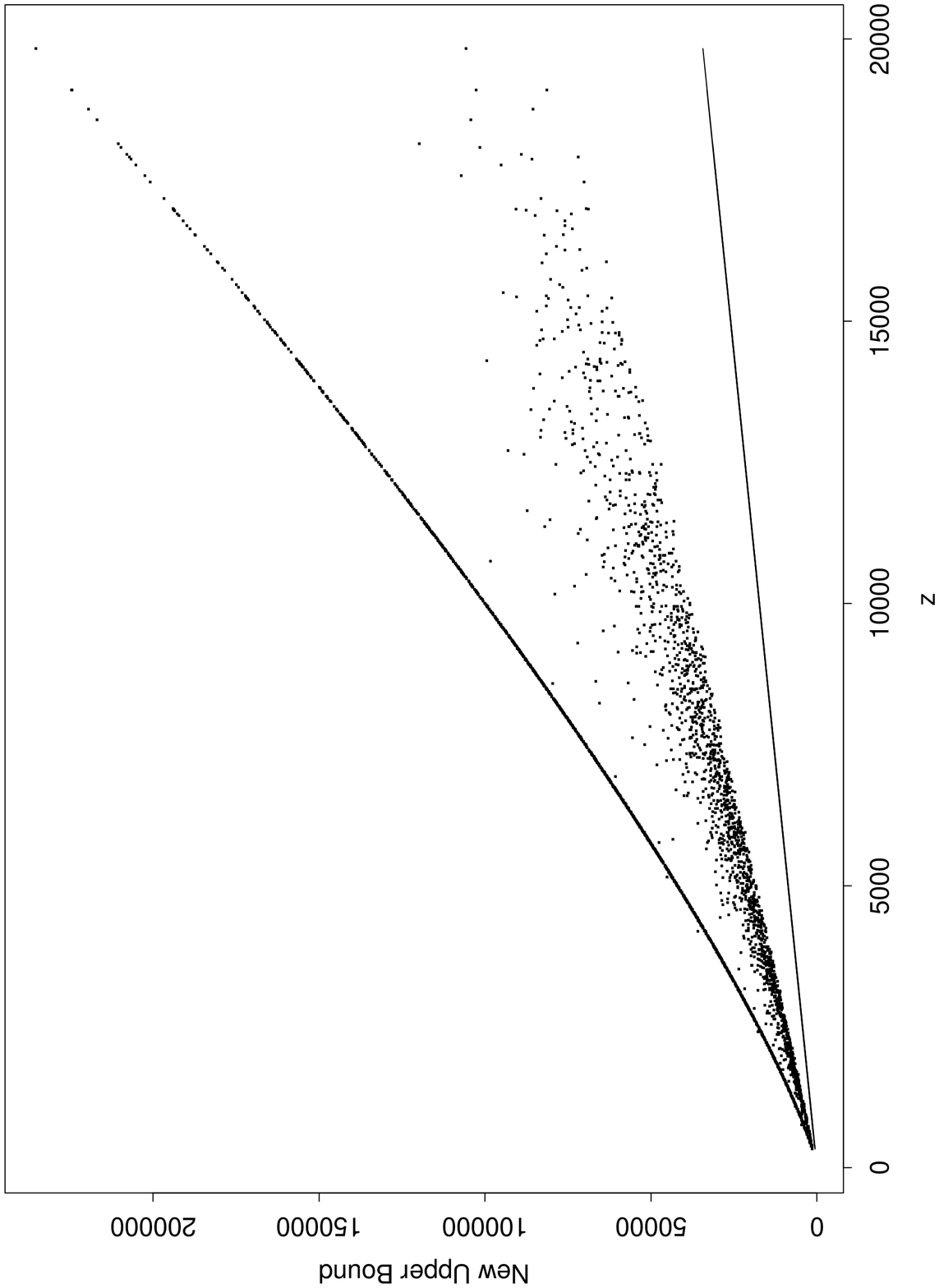}
\end{center}
\caption{The new upper bounds and the conjeture}\label{conjecturefig} 
\end{figure} 

In Figure \ref{conjecturefig} we plot the new upper bound as a function of $z$, and compare the points with 
Davison's lower bound (\ref{lowerbound}) and with the conjectured upper bound $z^{5/4}$. 
We see that most values of the new upper bound are smaller than $z^{5/4}$. This gives 
additional credence to the empirical conjecture in \cite{frobcomp}. 

\begin{figure}[htb]
\begin{center}
\includegraphics[totalheight=4in]{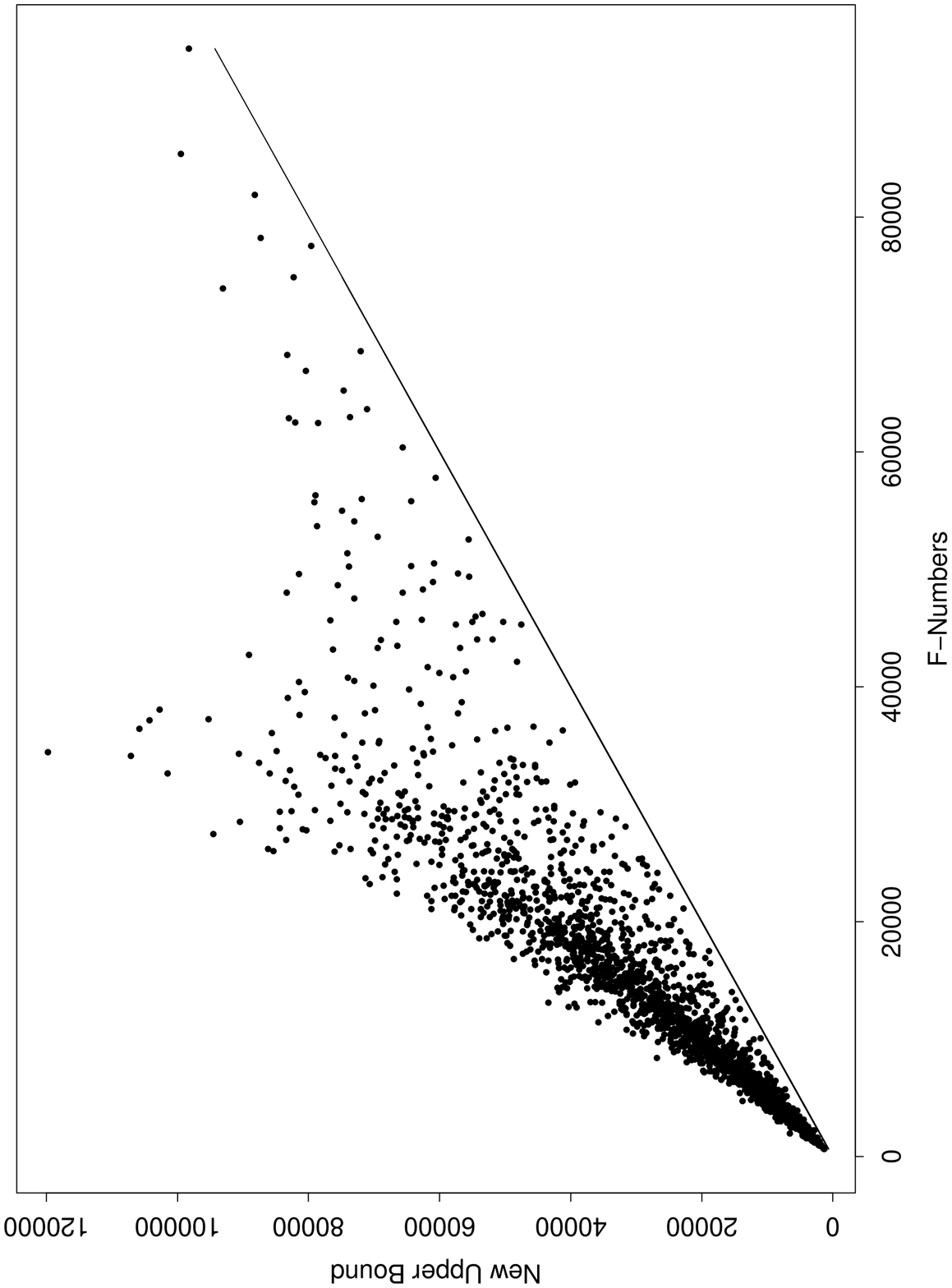}
\end{center}
\caption{The new upper bounds compared to the Frobenius numbers}\label{frobcomparisonfig} 
\end{figure} 

Finally, in Figure \ref{frobcomparisonfig} we plot the points of the new upper bound versus 
the true $f(a,b,c)$ values. We see that even for large values of $f(a,b,c)$ there are 
cases where the new upper bound yields close values. In 50\% of the cases the ratio of 
the new upper bound to the true $f(a,b,c)$ is greater than 2. 


\section{Final remarks} 

\begin{figure}[htb]
\begin{center}
\includegraphics[totalheight=4in]{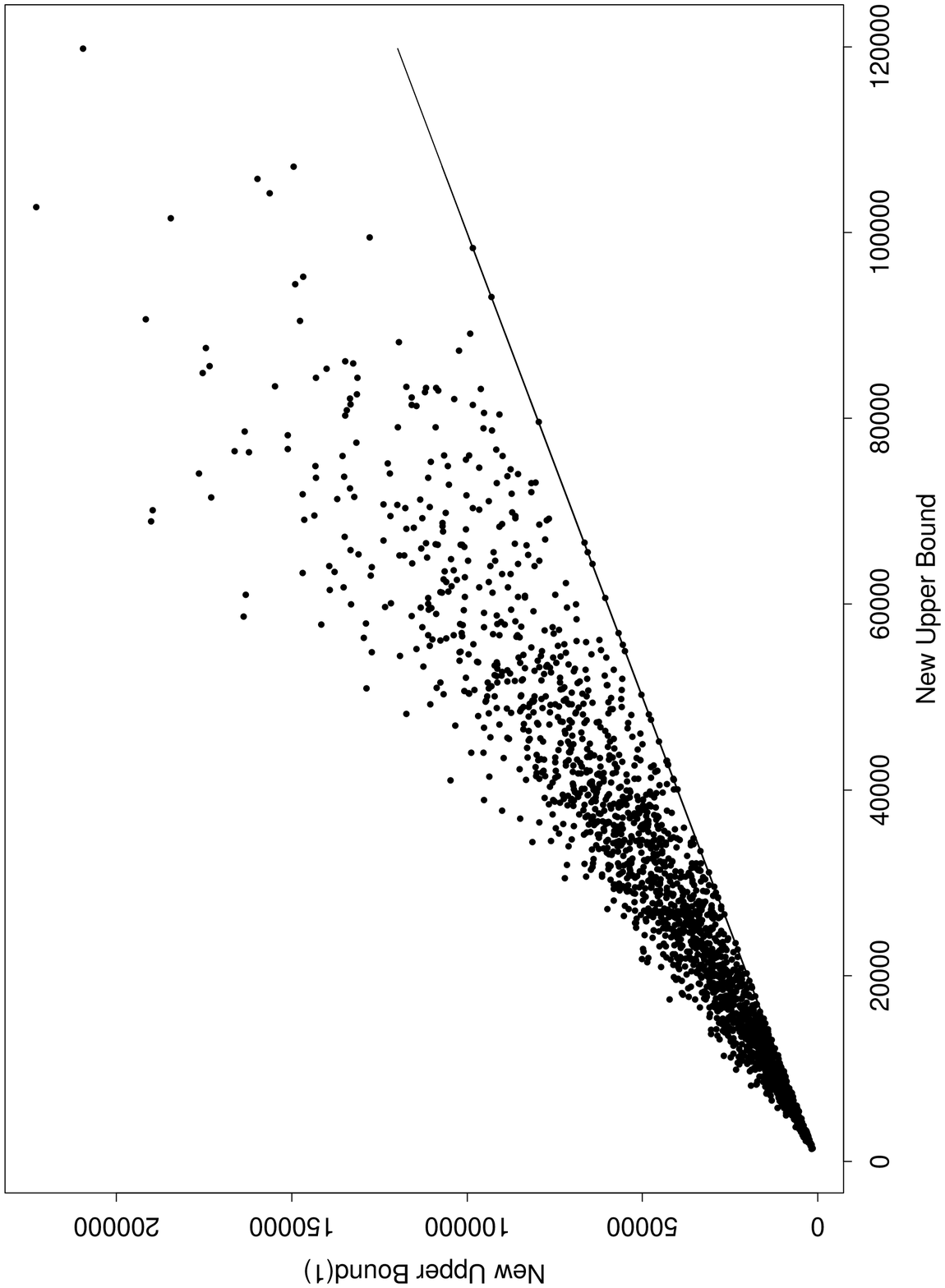}
\end{center}
\caption{The new upper bounds compared (1 iteration vs. 2 iterations)}\label{iterationcompfig} 
\end{figure}

As stated in the last section, we used ``only" two iterations ($N=2$) in our algorithm to compute bounds 
for $\sigma_t(a,b,c)$, which in turn lead to bounds on the Frobenius number $f(a,b,c)$. It is interesting 
to compare these values with the ones we get when just using the proposition, that is, one iteration ($N=1$). 
Figure \ref{iterationcompfig} illustrates this comparison. 

It is reasonable to expect even better results when one uses more than two iterations in the algorithm. 
However, we found that---at least for our range of variables---in the vast majority of cases 
the algorithm terminates prematurely after one or 
two iterations; accordingly, there is not much gain from increasing the number $N$ of iterations. 

The question which remains, even with a higher number of iterations, is how our bound can possibly be 
improved. The quality of the upper bound for $f(a,b,c)$ clearly depends only on the quality of the lower 
bound for $\sigma_t(a,b,c)$, and this bound is computed in our algorithm. There are three steps in the 
algorithm where we use bounds, namely, when adding[subtracting] $\Ql$[$\Qu$], respectively, in the second to last 
step where we use the Cauchy-Schwartz inequality, and in the last step where we use (\ref{sigmatrivbound}). 
Let us assume that we use enough iterations so that we only leave the REPEAT loop when $c_1 = 1$ or $c_2 = 1$; 
in this case we don't use the Cauchy-Schwartz inequality. The very last step does not entail a crude inequality, 
as we might be off maximally by one. This leaves the bounds $\Ql$ and $\Qu$ given by (\ref{Qlower}) and (\ref{Qupper}), 
respectively, and this is where we lose quality: here we might be off by a fractional factor of $c/c_1$, which then 
gets multiplied in the theorem, with possibly grave consequences. 

We close with a remark on the general Frobenius problem, that is, with an arbitrary number of arguments $a_1, \dots, a_d$. 
There are formulas analogous to (\ref{p_abc}) for $d>3$ \cite{bdr}; they involve higher-dimensional 
analogs of Dedekind-Rademacher sums. However, it is not clear how to compute them efficiently. 
According to \cite{barvinokalgorithm}, it is possible to compute these counting functions quickly; we 
just don't know an actual way to do so. It is our hope that the ideas in this paper can be extended 
to this general setting. 


\bibliographystyle{plain} 
\bibliography{bib} 

\setlength{\parskip}{0cm}

\end{document}